\documentclass{amsart}
\usepackage{amssymb}
\usepackage{amsmath}
\usepackage{amsfonts}

\newtheorem{theorem}{Theorem}
\theoremstyle{plain}

\newtheorem{corollary}{Corollary}

\newtheorem{definition}{Definition}

\newtheorem{lemma}{Lemma}
\newtheorem{notation}{Notation}

\newtheorem{proposition}{Proposition}
\newtheorem{remark}{Remark}

\numberwithin{equation}{section}
\input{tcilatex}

\begin{document}
\title[Random fixed point theorems]{Random fixed point theorems under mild
continuity assumptions}
\author{Monica Patriche}
\address{University of Bucharest\\
Departament of Mathematics, 14 Academiei Street, 010014, Bucharest}
\email{monica.patriche@yahoo.com}
\urladdr{http://fmi.unibuc.ro/ro/departamente/matematica/patriche\_monica/}
\dedicatory{}
\thanks{}

\begin{abstract}
In this paper, we study the existence of the random fixed points under mild
continuity assumptions. The main theorems consider the almost lower
semicontinuous operators defined on Fr\'{e}chet spaces and also operators
having properties weaker than lower semicontinuity. Our results either
extend or improve corresponding ones present in literature.
\end{abstract}

\maketitle

\textbf{\ }

\section{\textbf{INTRODUCTION}}

Fixed point theorems are a very powerful tool of the current mathematical
applications. These have been extended and generalized to study a wide class
of problems arising in mechanics, physics, economics and equilibrium theory,
engineering sciences, etc. New results concerning the deterministic or
random case were obtained, for instance, in \cite{fie1}, \cite{fie2}, \cite%
{pat1},\cite{pat2},\cite{pat3},\cite{pat4},\cite{pet1},\cite{pet},\cite{sha1}%
, \cite{sha2}, \cite{sha3}, \cite{sha4}.

The main aim of this work is to establish random fixed point theorems under
mild continuity assumptions. Results in this direction have been also
obtained, for example, in \cite{fie1}, \cite{fie2}, \cite{sha5}. Our
research enables us to improve some theorems obtained recently. We prove the
existence of the random fixed points for the lower semicontinuous operators
and for operators having properties weaker than lower semicontinuity. By
using the approximation method which is due to Ionescu-Tulcea (see \cite{ion}%
), we provide a new proof for the random Ky Fan fixed point theorem. We\
also introduce condition $\alpha $ which assures the existence of the random
fixed points$.$ The equivalence of condition $\alpha $ and condition ($%
\mathcal{P}$) (see \cite{fie1}) is proven and the assumptions over the
operator, which induce property $\alpha ,$ are found. Since the majority of
our results are obtained for operators defined on Fr\'{e}chet spaces, we
refer the reader to the new literature concerning this topic. For instance,
the authors worked also on Fr\'{e}chet spaces in \cite{aga1}, \cite{aga2},%
\cite{aga3}, \cite{don},\cite{don2},\cite{sha6}.

The rest of the paper is organized as follows. In the following section,
some notational and terminological conventions are given. We also present,
for the reader's convenience, some results on continuity and measurability
of the operators. Section 3 contains the main results. Section 4 presents
the conclusions of our research.

\section{\textbf{NOTATION AND DEFINITION}}

Throughout this paper, we shall use the following notation:

$2^{D}$ denotes the set of all non-empty subsets of the set $D$. If $%
D\subset Y$, where $Y$ is a topological space, cl$D$ denotes the closure of $%
D$. A \textit{paracompact space} is a Hausdorff topological space in which
every open cover admits an open locally finite refinement. Every metrizable
topological space is paracompact.

For the reader's convenience, we review a few basic definitions and results
from continuity and measurability of correspondences.

Let $Z$, $Y$ be topological spaces and $T:Z\rightarrow 2^{Y}$ be a
correspondence. $T$ is said to be \textit{upper semicontinuous} if for each $%
z\in Z$ and each open set $V$ in $Y$ with $T(z)\subset V$, there exists an
open neighborhood $U$ of $z$ in $Z$ such that $T(y)\subset V$ for each $y\in
U$. $T$ is said to be \textit{lower semicontinuous} if for each $z\in Z$ and
each open set $V$ in $Y$ with $T(z)\cap V\neq \emptyset $, there exists an
open neighborhood $U$ of $z$ in $Z$ such that $T(y)\cap V\neq \emptyset $
for each $y\in U$. $T:X\rightarrow 2^{Y}$ has open lower sections if $%
T^{-1}(y):=\{x\in X:y\in T(x)\}$ is open in $X$ for each $y\in Y.$ A
correspondence with open lower sections is lower semicontinuous.

Let $(X,d)$ be a metric space, $C$ be a \ nonempty subset of $X$ and $%
T:C\rightarrow 2^{X}$ be a correspondence.

We will use the following notations. We denote by $B(x,r)=\{y\in
C:d(y,x)<r\}.$ If $B_{0}$ is a subset of $X,$ then we will denote $%
B(B_{0},r)=\{y\in C:d(y,B_{0})<r\},$ where $d(y,B_{0})=\inf_{x\in B}d(y,x).$

We say that $T$ is \textit{hemicompact} if each sequence $\{x_{n}\}$ in $C$
has a convergent subsequence, whenever $d(x_{n}$ $T(x_{n}))\rightarrow 0$ as 
$n\rightarrow \infty .$\medskip

Let now $(\Omega ,\tciFourier $, $\mu )$ be a complete, finite measure
space, and $Y$ be a topological space. The correspondence $T:\Omega
\rightarrow 2^{Y}$ is said to have a \textit{measurable graph} if $G_{T}\in
\tciFourier \otimes \beta (Y)$, where $\beta (Y)$ denotes the Borel $\sigma $%
-algebra on $Y$ and $\otimes $ denotes the product $\sigma $-algebra. The
correspondence $T:\Omega \rightarrow 2^{Y}$ is said to be \textit{lower
measurable} if for every open subset $V$ of $Y$, the set $T^{-1}(V)=\{\omega
\in \Omega $ $:$ $T(\omega )\cap V\neq \emptyset $\} is an element of $%
\tciFourier $. This notion of measurability is also called in the literature 
\textit{weak measurability} or just \textit{measurability}, in comparison
with strong measurability: the correspondence $T:\Omega \rightarrow 2^{Y}$
is said to be \textit{strong measurable} if for every closed subset $V$ of $%
Y $, the set $\{\omega \in \Omega $ $:$ $T(\omega )\cap V\neq \emptyset $\}
is an element of $\tciFourier $. In the case when $X$ is separable, the
strong measurability coincides with the lower measurability.

Recall (see Debreu \cite{deb2}, p. 359) that if $T:\Omega \rightarrow 2^{Y}$
has a measurable graph, then $T$ is lower measurable. Furthermore, if $%
T(\cdot )$ is closed valued and lower measurable, then $T:\Omega \rightarrow
2^{Y}$ has a measurable graph.

A mapping $T:\Omega \times X\rightarrow Y$ is called a \textit{random
operator} if, for each $x\in X$, the mapping $T(\cdot ,x):\Omega \rightarrow
Y$ is measurable. Similarly, a correspondence $T:\Omega \times X\rightarrow
2^{Y}$ is also called a random operator if, for each $x\in X$, $T(\cdot
,x):\Omega \rightarrow 2^{Y}$ is measurable. A measurable mapping $\xi
:\Omega \rightarrow Y$ is called a \textit{measurable selection of the
operator} $T:\Omega \rightarrow 2^{Y}$ if $\xi (\omega )\in T(\omega )$ for
each $\omega \in \Omega $. A measurable mapping $\xi :\Omega \rightarrow Y$
is called a \textit{random fixed point} of the random operator $T:\Omega
\times X\rightarrow Y$ (or $T:\Omega \times X\rightarrow 2^{Y})$ if for
every $\omega \in \Omega $%
\c{}
$\xi (\omega )=T(\omega ,\xi (\omega ))$ $($or $\xi (\omega )=T(\omega ,\xi
(\omega ))$.

We will need the following measurable selection theorem in order to prove
our results.

\begin{proposition}
(Kuratowski-Ryll-Nardzewski Selection Theorem \cite{k}). A weakly measurable
correspondence with nonempty closed values from a measurable space into a
Polish space admits a measurable selector.
\end{proposition}

\section{MAIN\ RESULTS}

The aim of this section is to obtain new random fixed point theorems. The
main result concerns the operators which fulfill almost lower semicontinuity
or lower semicontinuity assumption. We extend the results from \cite{fie1}
in two aspects. Firstly, the assumptions of continuity are weakened and
secondly, the condition ($\mathcal{P}$)$,$ essential in the proof of the
existence of the random fixed points in the quoted paper is reformulated and
new assumptions which induce it are found.\medskip

\subsection{RANDOM FIXED POINT THEOREMS FOR OPERATORS WITH PROPERTIES WEAKER
THAN LOWER SEMICONTINUITY}

This subsection is mainly dedicated to establishing the random fixed point
theorems concerning the almost lower semicontinuous operators and the other
types of operators having properties weaker than lower semicontinuity. Our
results are new in literature. They can be compared with the ones stated in 
\cite{sha5}.

Firstly we recall the following proposition, which will be useful to prove
the main result of this subsection.

\begin{proposition}
(Theorem 3.4.in \cite{fie1}) Let $C$ be a closed separable subset of a
complete metric space $X$, and let $T:%
\Omega
\times C\rightarrow C(X)$ be a continuous hemicompact random operator. If,
for each $\omega \in 
\Omega
,$ the set $F(\omega ):=\{x\in C:x\in T(\omega ,x)\}\neq \emptyset $, then $%
T $ has a random fixed point.
\end{proposition}

Now, we are presenting the almost lower semicontinuous correspondences.

A correspondence $T:X\rightarrow 2^{Y}$ is \textit{almost lower
semicontinuous (a.l.s.c.)} \textit{at} $x$ (see \cite{deu}), if for any $%
\varepsilon >0$ there exists a neighborhood $U(x)$ of $x$ such that $%
\tbigcap\limits_{z\in U(x)}B(T(z);\varepsilon )\neq \emptyset .$

T is almost lower semicontinuous, if it is a.l.s.c. at each $x$.

If $\Omega $ is a nonempty set, we say that the operator $T:\Omega \times
X\rightarrow 2^{Y}$ is almost lower semicontinuous if for each $\omega \in
\Omega ,$ $T(\omega ,\cdot )$ is \ almost lower semicontinuous.

We note that when $Y$ is a normed linear space equipped with a norm $%
\left\Vert \cdot \right\Vert $, we may denote the $\varepsilon -$%
neighborhood of a subset $A$ of $Y$ by $B_{\varepsilon }(A):=\{y\in
Y:d(y,A)<\varepsilon \}$ where $d(y,A)=\inf \{\left\Vert x-y\right\Vert
:x\in A\}$.

In 1983, Deutsch and Kenderov \cite{deu} presented a remarkable
characterization of a.l.s.c. correspondences as follows.

\begin{lemma}
(see \cite{deu}) Let $X$ be paracompact, $Y$ a normat linear space, and $%
T:X\rightarrow 2^{Y}$ be a correspondence having convex images. Then T is
a.l.s.c: if and only if for each $\varepsilon >0$, $T$ admits a continuous $%
\varepsilon -$approximate selection f; that is, $f:X\rightarrow Y$ is a
continuous single-valued function such that $f(x)\in B(T(x);\varepsilon )$
for each $x\in X.$
\end{lemma}

The next theorem is the main result of this subsection. It states the
existence of the random fixed points for the almost lower semicontinuous
operators defined on Fr\'{e}chet spaces.

\begin{theorem}
Let $(\Omega ,\mathcal{F})$ be a measurable space, $C$ be a compact convex
separable subset of a Fr\'{e}chet space $X,$ and let $T:\Omega \times
C\rightarrow 2^{C}$ be a random operator. Suppose that for each $\omega \in
\Omega ,$ $T(\omega ,\cdot )$ is almost lower semicontinuous with non-empty
convex closed values and $(T(\omega ,\cdot ))^{-1}:\Omega \times
C\rightarrow 2^{C}$ is closed valued.
\end{theorem}

\textit{Then, }$T$\textit{\ has a random fixed point.}

\begin{proof}
Since for each $\omega \in \Omega ,$ $T(\omega ,\cdot )$ is almost lower
semicontinuous, according to Lemma 1, for each $n\in N,$\ there exists a
continuous function $f_{n}(\omega ,\cdot ):C\rightarrow C$\ such that $%
f_{n}(\omega ,x)\in B(T(\omega ,x);\frac{1}{n})$ for each\textit{\ }$x\in C.$
Let us define $T_{n}:\Omega \times C\rightarrow 2^{C}$ by $T_{n}(\omega
,x)=B(T(\omega ,x);\frac{1}{n})$ for each $(\omega ,x)\in \Omega \times C.$

According to Tihonov fixed point theorem, for each $n\in N,$ there exists $%
x_{n}\in C$ such that $x_{n}=f_{n}(\omega ,x_{n})$ and then, $x_{n}\in
B(T(\omega ,x_{n});\frac{1}{n}).$

$C$ is compact, then $f_{n}$ is hemicompact for each $n\in \mathbb{N}$.
According to Proposition 2, for each $n\in \mathbb{N},$ $f_{n}$ has a random
fixed point and then, $T_{n}$ has a random fixed point $\xi _{n},$ that is $%
\xi _{n}:\Omega \rightarrow C$ is measurable and $\xi _{n}(\omega )\in
T_{n}(\omega ,\xi _{n}(\omega ))$ for $n\in N$.

Let $\omega \in \Omega $ be fixed. Then, $d(\xi _{n}(\omega ),T(\omega ,\xi
_{n}(\omega ))\rightarrow 0$ when $n\rightarrow \infty $ and since $C$ is
compact, $\{\xi _{n}(\omega )\}$ has a convergent subsequence $\{\xi
_{n_{k}}(\omega )\}.$ Let $\xi _{0}(\omega )=\lim_{n\rightarrow \infty }\xi
_{n_{k}}(\omega ).$ It follows that $\xi _{0}:\Omega \rightarrow C$ is
measurable and for each $\omega \in \Omega ,$ $d(\xi _{0}(\omega ),T(\omega
,\xi _{n_{k}}(\omega ))\rightarrow 0$ when $n_{k}\rightarrow \infty .$

Let us assume that there exists $\omega \in \Omega $ such that $\xi
_{0}(\omega )\notin T(\omega ,\xi _{0}(\omega )).$ Since $\{\xi _{0}(\omega
)\}\cap (T(\omega ,\cdot ))^{-1}(\xi _{0}(\omega ))=\emptyset $ and $X$ is a
regular space, there exists $r_{1}>0$ such that $B(\xi _{0}(\omega
),r_{1})\cap (T(\omega ,\cdot ))^{-1}(\xi _{0}(\omega ))=\emptyset $.
Consequently, for each $z\in B(\xi _{0}(\omega ),r_{1})$ we have that $%
z\notin (T(\omega ,\cdot ))^{-1}(\xi _{0}(\omega ))$ which is equivalent
with $\xi _{0}(\omega )\notin T(\omega ,z)$ or $\{\xi _{0}(\omega )\}\cap
T(\omega ,z)=\emptyset $. The closedness of each $T(\omega ,z)$ and the
regularity of $X$ implies the existence of a real number $r_{2}>0$ such that 
$B(\xi _{0}(\omega ),r_{2})\cap T(\omega ,z)=\emptyset $ for each $z\in
B(\xi _{0}(\omega ),r_{1}),$ which implies $\xi _{0}(\omega )\notin
B((\omega ,z);r_{2})$ for each $z\in B(\xi _{0}(\omega ),r_{1}).$ Let $%
r=\min \{r_{1},r_{2}\}.$ Hence, $\xi _{0}(\omega )\notin B((\omega ,z);r)$
for each $z\in B(\xi _{0}(\omega ),r)$ and then, there exists $N^{\ast }\in 
\mathbb{N}$ such that for each $n_{k}>N^{\ast },$ $\xi _{0}(\omega )\notin
B((\omega ,\xi _{n_{k}}(\omega ));r)$ which contradicts $d(\xi _{0}(\omega
),T(\omega ,\xi _{n_{k}}(\omega ))\rightarrow 0$. It follows that our
assumption is false.

Hence for each $\omega \in \Omega ,$ $\xi _{0}(\omega )\in T(\omega ,\xi
_{0}(\omega ))$, where $\xi _{0}:\Omega \rightarrow C$ is measurable. We
conclude that $T$ has a random fixed point.
\end{proof}

Related to the almost lower semicontinuous correspondences, there are
correspondences with local intersection property and sub-lower
semicontinuous correspondences, which differ very slightly from the first
ones. We will obtain some results related to Theorem 1 in these cases and
will introduce the reader in the topic we brought into discussion by
presenting firstly the definitions.\medskip

A corollary can be proved for the operators which satisfy the local
intersection property, which is defined below.

Let $X$, $Y$ be topological spaces. The correspondence $T:X\rightarrow 2^{Y}$
has the \textit{local} \textit{intersection property (see }\cite{wu}) if $%
x\in X$ with $T(x)\neq \emptyset $ implies the existence of an open
neighborhood $V(x)$ of $x$ such that $\cap _{z\in V(x)}T(z)\neq \emptyset .$

If $\Omega $ is a nonempty set, we say that the operator $T:\Omega \times
X\rightarrow 2^{Y}$ has local intersection property if for each $\omega \in
\Omega ,$ $T(\omega ,\cdot )$ has the local intersection property.\textit{%
\medskip }

We establish the next theorem, which concerns operators having the local
intersection property. Its proof relies on the following lemma.

\begin{lemma}
(Wu and Shen \cite{wu}). \textit{Let }$X$\textit{\ be a non-empty
paracompact subset of a Hausdorff topological space }$E$\textit{\ and }$Y$%
\textit{\ be a non-empty subset of a Hausdorff topological vector space. Let 
}$S$\textit{,}$T:X\rightarrow 2^{Y}$\textit{\ be correspondences which
verify:}
\end{lemma}

1)\textit{\ for each }$x\in X,$\textit{\ }$S(x)\neq \emptyset $\textit{\ and 
}$coS(x)\subset T\left( x\right) ;$

2)\textit{\ }$S$\textit{\ has local intersection property.}

\textit{Then, }$T$\textit{\ has a continuous selection.\medskip }

The corresponding random fixed point theorem in case of local intersection
property is stated below.

\begin{theorem}
Let $(\Omega ,\mathcal{F})$ be a measurable space, $C$ be a compact convex
separable subset of a Fr\'{e}chet space $X,$ and let $T:\Omega \times
C\rightarrow 2^{C}$ be a random operator. Suppose, for each $\omega \in
\Omega ,$ that $T(\omega ,\cdot )$ has the local intersection property and
convex values.
\end{theorem}

\textit{Then, }$T$\textit{\ has a random fixed point. }\medskip

\begin{proof}
Since for each $\omega \in \Omega ,$ $T(\omega ,\cdot )$ has local
intersection property, according to Lemma 2,\ there exists a continuous
function $f(\omega ,\cdot ):C\rightarrow C$\ such that for each\textit{\ }$%
x\in C,$\textit{\ }$f(\omega ,x)\in T(\omega ,x).$

According to Tihonov's fixed point Theorem, there exists $x\in C$ such that $%
x=f(\omega ,x)$ and then, $x\in T(\omega ,x).$

$C$ is compact, then $f$ is hemicompact. According to Proposition 2, $f$ has
a random fixed point and then, $T$ has a random fixed point $\xi ,$ that is $%
\xi :\Omega \rightarrow C$ is measurable and $\xi (\omega )\in T(\omega ,\xi
(\omega ))$.
\end{proof}

The sub-lower semicontinuous correspondences were defined by Zheng in \cite%
{142}.

Let $X$ be a topological space and $Y$ be a topological vector space, a
correspondence $T:X\rightarrow 2^{Y}$ is called \textit{sub-lower
semicontinuous} \cite{142} if for each $x\in X$ and for each neighborhood $V$
of $0$ in $Y,$ there exists $z\in T(x)$ and a neighborhood $U(x)$ of $x$ in $%
X$ such that for each $y\in U(x),$ $z\in T(y)+V.$

If $\Omega $ is a nonempty set, we say that the operator $T:\Omega \times
X\rightarrow 2^{Y}$ is sub-lower semicontinuous if for each $\omega \in
\Omega ,$ $T(\omega ,\cdot )$ is sub-lower semicontinuous.\medskip

Zheng proved in \cite{142} a continuous selection result for the sub-lower
semicontinuous correspondences, which can be used in order to obtain Theorem
3. Here is his result.

\begin{lemma}
\cite{142} \textit{Let }$X$\textit{\ be a paracompact topological space, }$Y$%
\textit{\ be a locally convex topological linear space and let }$%
T:X\rightarrow 2^{Y}$\textit{\ be a correspondence with convex values. Then, 
}$T$\textit{\ is sub-lower semicontinuous\ if and only if for each
neighborhood }$V$\textit{\ of }$0$\textit{\ in }$Y,$\textit{\ there exists a
continuous function }$f:X\rightarrow Y$\textit{\ such that for each }$x\in
X, $\textit{\ }$f(x)\in T(x)+V.$
\end{lemma}

The random fixed point existence for the sub-lower semicontinuous random
operators is stated below.

\begin{theorem}
Let $(\Omega ,\mathcal{F})$ be a measurable space, $C$ be a compact convex
separable subset of a Fr\'{e}chet space $X,$ and let $T:\Omega \times
C\rightarrow 2^{C}$ be a random operator. Suppose, for each $\omega \in
\Omega ,$ that $T(\omega ,\cdot )$ is sub-lower semicontinuous with
non-empty convex closed values and $(T(\omega ,\cdot ))^{-1}:\Omega \times
C\rightarrow 2^{C}$ is closed valued.
\end{theorem}

\textit{Then, }$T$\textit{\ has a random fixed point.\medskip }

The proof of Theorem 3 is similar to the one of Theorem 1, but it relies on
Lemma 3.\medskip 

In \cite{ans}, Ansari and Yao proved a fixed point theorem for transfer
open-valued correspondences. Their proof is based on a continuous selection
theorem which they construct. We present below their result.\medskip

Let $X$ and $Y$ be two topological vector spaces. The correspondence $%
T:X\rightarrow 2^{Y}$ is said to be \textit{transfer open-valued} (see \cite%
{ans}) if for any $x\in X$, $y\in T(x),$ there exists an $z\in X$ such that $%
y\in $int$_{Y}T(z)$.

If $\Omega $ is a nonempty set, we say that the operator $T:\Omega \times
X\rightarrow 2^{Y}$ is transfer open valued if for each $\omega \in \Omega ,$
$T(\omega ,\cdot )$ is transfer open valued.\medskip

The proof of the next lemma is included in the proof of Theorem 1 in \cite%
{ans}.

\begin{lemma}
Let $K$ be a non-empty subset of a Hausdorff topological vector space $E$,
and let $T:K\rightarrow 2^{K}$ be a correspondence with non empty convex
values. If $K=\tbigcup \{$int$_{K}T^{-1}(y):y\in K\}$ (or $T^{-1}$ is
transfer open valued), then $T$ has a continuous selection.
\end{lemma}

We obtain Theorem 4, by using Lemma 4 and the same argument as in the proof
of Theorem 2.

\begin{theorem}
Let $(\Omega ,\mathcal{F})$ be a measurable space, $C$ be a compact convex
separable subset of a Fr\'{e}chet space $X,$ and let $T:\Omega \times
C\rightarrow 2^{C}$ be a random operator with non-empty convex values.
Suppose that $K=\tbigcup \{$int$_{K}(T(\omega ,\cdot ))^{-1}(y):y\in K\}$ or
for each $\omega \in \Omega ,$ $(T(\omega ,\cdot ))^{-1}$ is transfer open
valued.
\end{theorem}

\textit{Then, }$T$\textit{\ has a random fixed point.\medskip }

We also state the following result.

\begin{theorem}
Let $(\Omega ,\mathcal{F})$ be a measurable space, $C$ be a compact convex
separable subset of a Fr\'{e}chet space $X,$ and let $T:\Omega \times
C\rightarrow 2^{C}$ be a random operator with non-empty closed convex
values, such that for each $\omega \in \Omega ,$ $(T(\omega ,\cdot
))^{-1}:\Omega \times C\rightarrow 2^{C}$ is closed valued. Suppose that for
each open neighborhood $V$ of the origin and for each $\omega \in \Omega ,$
the correspondence $(S^{V,\omega })^{-1}:C\rightarrow 2^{C}$ is transfer
open valued, where $S^{V,\omega }(x)=(T(\omega ,x)+V)\cap C$ for each $x\in
C $.
\end{theorem}

\textit{Then, }$T$\textit{\ has a random fixed point.\medskip }

\begin{proof}
By using Lemma 4, we prove that for each $n\in \mathbb{N}$ and for each $%
\omega \in \Omega $ there exists a continuous function $f_{n}:C\rightarrow C,
$ $f_{n}(\omega ,x)\in B(T(\omega ,x);\frac{1}{n})\cap C$ for each $x\in C.$
The proof in similar to the proof of Theorem 1.
\end{proof}

\subsection{RANDOM FIXED POINT THEOREMS FOR LOWER SEMICONTINUOUS OPERATORS}

This subsection is designed to extending the results established in \cite%
{fie1} by considering lower semicontinuous operators defined on Fr\'{e}chet
spaces.

Condition ($\mathcal{P}$), firstly introduced by Petryshyn \cite{p1} in
order to prove the existence of fixed points for single-valued operators was
extended to multivalued operators. We provide here the definition for the
last case.

Let $(X,d)$ be a metric space, $C$ be a \ nonempty closed subset of $X$ and $%
T:C\rightarrow 2^{X}$ be a correspondence. $T$ is said to satisfy \textit{%
condition} ($\mathcal{P}$) (see \cite{fie1}) if, for every closed ball $B$
of $C$ with radius $r\geq 0$ and any sequence $\{x_{n}\}$ in $C$ for which $%
d(x_{n},B)\rightarrow 0$ and $d(x_{n},T(x_{n}))\rightarrow 0$ as $%
n\rightarrow \infty ,$ there exists $x_{0}\in B$ such that $x_{0}\in
T(x_{0}).$ If $\Omega $ is any nonempty set, we say that the operator $%
T:\Omega \times C\rightarrow 2^{X}$ satisfies condition ($\mathcal{P}$) if
for each $\omega \in \Omega ,$ the correspondence $T(\omega ,\cdot
):C\rightarrow 2^{X}$ satisfies condition ($\mathcal{P}$).\medskip

We also present the main result in \cite{fie1}, concerning operators which
satisfy condition ($\mathcal{P}$). We will extend further this theorem.

\begin{proposition}
(Theorem 3.2 in \cite{fie1}) Let $C$ be a closed separable subset of a
complete metric space and let $T:\Omega \times C\rightarrow 2^{X}$ be a
lower semi-continuous random operator, which enjoys condition ($\mathcal{P}$%
). Suppose, for each $\omega \in \Omega $ that the set
\end{proposition}

$F(\omega ):=\{x\in C:x\in T(\omega ,x)\}\neq \emptyset .$

\textit{Then }$T$\textit{\ has a random fixed point theorem.\medskip }

Fierro, Mart\i nez and Morales showed in \cite{fie1} (the proof of Theorem
3.1) that if $T$ is lower semicontinuous, then it satisfies the following
condition, which we call (*), condition necessary to prove the existence of
random fixed points.

(*): \ \ \ if $B_{0}$ is a closed ball and there exists $x_{0}\in B_{0}$
such that $x_{0}\in T(x_{0}),$

\ \ \ \ \ \ \ \ then, there exists a convergent sequence $x_{n}\in
B(B_{0};r) $ such that

$\ \ \ \ \ \ \ x_{n}\rightarrow x_{0}$ when $n\rightarrow \infty $ and $%
d(x_{n},T(x_{n}))<\frac{1}{n}.$\textit{\medskip }

Our work will consider simpler assumptions which imply condition (*).\medskip

Firstly we introduce condition $\alpha $. The following lemmata establish
the connection with condition ($\mathcal{P}$), widely used in old results
concerning random fixed points and we state simpler assumptions which
assures the fulfillment of condition $\alpha .$ These results will be used
in order to establish new random fixed point theorems, which will be
compared with the existent ones in literature.

\begin{definition}
Let $(X,d)$ be a metric space and $C$ a \ nonempty subset of $X.$
\end{definition}

\textit{1) We say that the correspondence }$T:C\rightarrow 2^{X}$\textit{\
satisfies condition }$\alpha $\textit{\ if }$x_{0}\notin T(x_{0})$\textit{\
implies the existence of a real }$r>0$\textit{\ such that }$x_{0}\notin
B(T(x);r)\cap B(x,r)$\textit{\ for each} $x\in B(x_{0},r).$

2) \textit{If }$\Omega $\textit{\ is a nonempty set, we say that the
operator }$T:\Omega \times C\rightarrow 2^{X}$\textit{\ satisfies condition }%
$\alpha $\textit{\ if, for each }$\omega \in \Omega ,$\textit{\ the
correspondence }$T(\omega ,\cdot ):C\rightarrow 2^{X}$\textit{\ satisfies
condition }$\alpha .\medskip $

The equivalence between condition ($\mathcal{P}$) and condition $\alpha $
for correspondences defined on locally compact subsets of a complete metric
space is showed in Lemma 5.

\begin{lemma}
Let $(X,d)$ be a complete metric space, $C$ be a nonempty closed separable
subset of $X$ and $T:C\rightarrow 2^{X}$ be a correspondence.
\end{lemma}

\textit{If }$T$\textit{\ satisfies condition (}$\mathcal{P}$\textit{), then,
it satisfies condition }$\alpha .$ \textit{If }$C$\textit{\ is locally
compact and }$T$\textit{\ satisfies condition (}$\alpha $\textit{), then, }$%
T $\textit{\ satisfies condition (}$\mathcal{P}$\textit{)}$.$

\begin{proof}
Let us consider a correspondence $T$ which satisfies condition ($\mathcal{P}$%
). We will prove that $T$ satisfies condition $\alpha .$ Let $x_{0}\in C$ be
a point such that $x_{0}\notin T(x_{0}).$ Let us assume, by contrary, that
for each $r>0,$ there exists $x_{r}\in B(x_{0},r)$ such that\textit{\ }$%
x_{0}\in B(T(x_{r});r)\cap B(x_{r},r).$ Therefore, for each natural number $%
n>0$\ there exists $x_{n}\in B(x_{0},\frac{1}{n})$ such that\textit{\ }$%
x_{0}\in B(T(x_{n});\frac{1}{n})\cap B(x_{n},\frac{1}{n}).$\textit{\ }%
Consequently,\textit{\ }we found a sequence $\{x_{n}\}$ in $C$ with the
property that $d(x_{n},x_{0})\rightarrow 0$ and $d(x_{n},T(x_{n}))%
\rightarrow 0.$ Since $T$ satisfies condition ($\mathcal{P}$) for $%
B_{0}=\{x_{0}\},$ it follows that $x_{0}\in T(x_{0}),$ which contradicts $%
x_{0}\notin T(x_{0}).$ This means that our assumption is false, and then,
condition ($\mathcal{P}$) implies condition $\alpha .$

For the inverse implication, let us consider a closed ball $B_{0}$ of $C$
with radius $r\geq 0$ and a sequence $\{x_{n}\}$ in $C$ for which $%
d(x_{n},B_{0})\rightarrow 0$ and $d(x_{n},T(x_{n}))\rightarrow 0$ as $%
n\rightarrow \infty .$ We use the sequentially compactness of the set $G$ we
will define in the following way. Let us firstly denote for each $n\in N,$ $%
r_{n}=d(x_{n},B_{0}).$ According to the hypotheses, the sequence $\{r_{n}\}$
is convergent and $\lim_{n\rightarrow \infty }r_{n}=0.$ Let us consider $%
\varepsilon >0$ and then, there exists $N(\varepsilon )\in \mathbb{N}$ such
that $r_{n}<\varepsilon $ for each $n>N(\varepsilon ).$ Since $X$ is locally
compact, the small closed balls are compact and consequently, the set $G$
defined as $G=$cl$((B_{0})_{\varepsilon }\cap C)$ is compact. Therefore, the
sequence $\{x_{n}\}\subset G$ has a convergent subsequence $\{x_{n_{k}}\}$
and, without loss of generality, we can assume that the sequence $\{x_{n}\}$
is convergent. It remains to show the following assertion:

If $T:C\rightarrow 2^{X}$ satisfies condition $\alpha $, and if there exists
a closed ball $B_{0}$ and a convergent sequence $\{x_{n}\}\rightarrow x_{0}$
such that $d(x_{n},B_{0})\rightarrow 0$ and $d(x_{n},T(x_{n}))\rightarrow 0$
then, $x_{0}\in T(x_{0})$ and $x_{0}\in B_{0}.$

In order to prove this, we note that the closedness of $B_{0}$ implies $%
x_{0}\in B_{0}.$ Let us assume, by contrary, that $x_{0}\notin T(x_{0}).$
Then, according to condition $\alpha ,$ there exists $r>0$ such that $%
x_{0}\notin B(T(x);r)\cap B(x,r)$ for each $x\in B(x_{0},r).$ The
convergence of $\{x_{n}\}$ to $x_{0}$ implies the existence of a natural
number $N(r)\in \mathbb{N}$ such that $x_{n}\in B(x_{0},r)$ for each $%
n>N(r). $ Consequently, $x_{0}\notin B(T(x_{n});r)\cap B(x_{n},r)$ for each $%
n>N(r).$ Since for each $n>N(r),$ $x_{0}\in B(x_{n},r),$ it follows that if $%
n>N(r),$ $x_{0}\notin B(T(x_{n});r),$ that is $d(x_{0},T(x_{n}))>r.$ This
fact contradicts $d(x_{0},T(x_{n}))\rightarrow 0$ when $n\rightarrow \infty
, $ which is true from the hypotheses and because $x_{0}=\lim_{n\rightarrow
\infty }x_{n}.$ This means that our assumption is false, and it remains that 
$x_{0}\in T(x_{0}).$ We proved that if $T$ satisfies condition $\alpha ,$
then $T$ satisfies condition ($P$).
\end{proof}

\begin{lemma}
Let $(X,d)$ be a metric space, $C$ be a nonempty closed subset of $X$ and $%
T:C\rightarrow 2^{X}$ be a correspondence.
\end{lemma}

\textit{If }$T$ \textit{and} $T^{-1}$\textit{\ have closed values, then, }$T$%
\textit{\ satisfies condition }$\alpha $. \textit{If, in addition, }$C$%
\textit{\ is locally compact, then }$T$\textit{\ satisfies condition (}$P$%
\textit{))}$.\medskip $

\begin{proof}
Let us consider $x_{0}\in C$ such that $x_{0}\notin T(x_{0}).$ Since $%
\{x_{0}\}\cap T^{-1}(x_{0})=\emptyset $ and $X$ is a regular space, there
exists $r_{1}>0$ such that $B(x_{0},r_{1})\cap
B(T^{-1}(x_{0});r_{1}))=\emptyset ,$ and then, $B(x_{0},r_{1})\cap
T^{-1}(x_{0})=\emptyset $. Consequently, for each $x\in B(x_{0},r_{1})$ we
have that $x\notin T^{-1}(x_{0})$ which is equivalent with $x_{0}\notin T(x)$
or $\{x_{0}\}\cap T(x)=\emptyset $. The closedeness of each $T(x)$ and the
regularity of $X$ implies the existence of a real number $r_{2}>0$ such that 
$B(x_{0},r_{2})\cap T(x)=\emptyset $ for each $x\in B(x_{0},r_{1}),$ which
implies $x_{0}\notin B(T(x);r_{2})$ for each $x\in B(x_{0},r_{1}).$ Let $%
r=\min \{r_{1},r_{2}\}.$ Hence, $x_{0}\notin B(T(x);r)$ for each $x\in
B(x_{0},r)$ and thus condition $\alpha $ is fulfilled.
\end{proof}

The following theorem states the existence of the random fixed points for
lower semicontinuous correspondences defined on locally compact complete
metric spaces.

\begin{theorem}
Let $(\Omega ,\mathcal{F})$ be a measurable space, $C$ be a closed separable
subset of a locally compact complete metric space and let $T:\Omega \times
C\rightarrow 2^{X}$ be a lower semicontinuous random operator with closed
values. Suppose that, for each $\omega \in \Omega ,$ $(T(\omega ,\cdot
))^{-1}:X\rightarrow 2^{C}$ is closed valued and the set
\end{theorem}

$F(\omega ):=\{x\in C:x\in T(\omega ,x)\}\neq \emptyset .$

\textit{Then }$T$\textit{\ has a random fixed point theorem.\medskip }

\begin{proof}
Since $C$ is locally compact and for each $\omega \in \Omega ,$ $T(\omega
.\cdot ):C\rightarrow 2^{X}$ and $(T(\omega .\cdot ))^{-1}:X\rightarrow
2^{C} $ have closed values, by applying Lemma 6, we obtain that $T$
satisfies condition ($\mathcal{P}$). All the assumptions of Proposition 3
are fulfilled, then, $T$\textit{\ }has a random fixed point.
\end{proof}

The existence of the random fixed points remains valid if for each $\omega
\in \Omega ,$ $(T(\omega .\cdot ))^{-1}:X\rightarrow 2^{C}$ is lower
semicontinuous.

\begin{theorem}
Let $(\Omega ,\mathcal{F})$ be a measurable space, $C$ be a closed separable
subset of a locally compact complete metric space and let $T:\Omega \times
C\rightarrow 2^{X}$ be a random operator with closed values. Suppose that
for each $\omega \in \Omega ,$ $(T(\omega .\cdot ))^{-1}:X\rightarrow 2^{C}$
is lower semicontinuous and closed valued and the set $F(\omega ):=\{x\in
C:x\in T(\omega ,x)\}\neq \emptyset .$
\end{theorem}

\textit{Then }$T$\textit{\ has a random fixed point theorem.\medskip }

\begin{proof}
According to Proposition 3, there exists a measurable mapping $\xi :\Omega
\rightarrow C$ such that\textit{\ }for each $\omega \in \Omega $%
\c{}
$\xi (\omega )\in (T(\omega ,\cdot ))^{-1}(\xi (\omega )),$ that is, for
each $\omega \in \Omega $%
\c{}
$\xi (\omega )\in T(\omega ,\xi (\omega )).$ Therefore, we obtained a random
fixed point for $T.$
\end{proof}

The next result, due to Michael, is very important in the theory of
continuous selections.

\begin{lemma}
(Michael \cite{mic}). \textit{Let }$X$\textit{\ be a }$T_{1}$\textit{,
paracompact space. If }$Y$\textit{\ is a Banach space, then each lower
semicontinuous convex closed valued correspondence }$T:X\rightarrow 2^{Y}$%
\textit{\ admits a continuous selection.}
\end{lemma}

By using the above lemma, we establish Theorem 8.

\begin{theorem}
Let $(\Omega ,\mathcal{F})$ be a measurable space, $C$ be a compact convex
separable subset of a compact Fr\'{e}chet space and let $T:\Omega \times
C\rightarrow 2^{C}$ be a lower semicontinuous random operator with closed
convex values, such that for each $\omega \in \Omega ,$ $(T(\omega .\cdot
))^{-1}:C\rightarrow 2^{C}$ is closed valued.
\end{theorem}

\textit{Then }$T$\textit{\ has a random fixed point theorem.\medskip }

\begin{proof}
Since $C$ is compact and for each $\omega \in \Omega ,$ $T(\omega .\cdot
):C\rightarrow 2^{C}$ and $(T(\omega .\cdot ))^{-1}:C\rightarrow 2^{C}$ have
closed values, by applying Lemma 6, we obtain that $T$ satisfies condition ($%
\mathcal{P}$).

According to Michael's selection theorem, there exists a continuous function 
$f:C\rightarrow C$ such that $f$ is continuous and $f(x)\in T(\omega ,x).$
According to Tihonov fixed point theorem, there exists deterministic fixed
points for $f,$ then, the set $F(\omega ):=\{x\in C:x\in T(\omega ,x)\}\neq
\emptyset .$ All the assumptions of Theorem 6 are fulfilled, then, $T$%
\textit{\ }has a random fixed point.
\end{proof}

We will provide another condition which implies (*).

\begin{definition}
Let $(X,d)$ be a complete metric space and $C$ a \ nonempty closed subset of 
$X.$
\end{definition}

\textit{1) We say that the correspondence }$T:C\rightarrow 2^{X}$\textit{\
satisfies condition }$\beta $ \textit{if }$x_{0}\in T(x_{0})$ \textit{%
implies that for each open neighborhood of the origin in }$X$\textit{, there
exists an open neighborhood }$U(x_{0})$ \textit{of }$x_{0}$\textit{\ such
that }$x_{0}\in T(x)+V$\textit{\ for each }$x\in U(x_{0}).$\textit{\ }

2) \textit{If }$\Omega $\textit{\ is a nonempty set, we say that the
operator }$T:\Omega \times C\rightarrow 2^{X}$\textit{\ satisfies condition }%
$\beta $\textit{\ if, for each }$\omega \in \Omega ,$\textit{\ the
correspondence }$T(\omega ,\cdot ):C\rightarrow 2^{X}$\textit{\ satisfies
condition }$\beta .\medskip $

\begin{remark}
If $T$ satisfies condition $\beta ,$ then, $T$ is sub-lower semicontinuous 
\cite{142} in its fixed points.
\end{remark}

\begin{lemma}
Let $(X,d)$ be a complete metric space and $C$ a nonempty closed subset of $%
X.$ \textit{If the correspondence }$T:C\rightarrow 2^{X}$ satisfies
condition $\beta $, then, $T$ satisfies (*).
\end{lemma}

\begin{proof}
Let $B_{0}$ be a closed ball and $x_{0}\in B_{0}$ such that $x_{0}\in
T(\omega ,x_{0}).$ According to condition $\beta $, for each $n\in \mathbb{N}
$, $x_{0}\in B(T(x),\frac{1}{n})$ for each $x\in B(x_{0},\frac{1}{n})\cap
U(x_{0}).$ We have that $B(x_{0},\frac{1}{n})\cap U(x_{0})\cap B(T(x),\frac{1%
}{n})\neq \emptyset $ for each $x\in B(x_{0},\frac{1}{n})\cap U(x_{0}).$
Since $B(x_{0},\frac{1}{n})\cap U(x_{0})$ and $B(T(x),\frac{1}{n})$ are open
sets, $B(x_{0},\frac{1}{n})\cap U(x_{0})\cap B(T(x),\frac{1}{n})\neq
\{x_{0}\}.$ Therefore, for each $n\in \mathbb{N}$, we can choose $x_{n}\in
B(x_{0},\frac{1}{n}),$ $x_{n}\neq x_{0}$ and then, $x_{n}\in B(B_{0}$,$\frac{%
1}{n}),$ $x_{n}\rightarrow x_{0}$ when $n\rightarrow \infty $ and $%
d(x_{n},T(x_{n}))<\frac{1}{n},$ that is $T$ satisfies (*).
\end{proof}

We obtain the following random fixed point theorem.

\begin{theorem}
Let $(\Omega ,\mathcal{F})$ be a measurable space, $C$ be a closed separable
subset of a complete metric space and let $T:\Omega \times C\rightarrow
2^{X} $ be a random operator which enjoys conditions ($\mathcal{P}$) and $%
\beta $. Suppose, for each $\omega \in \Omega $ that the set
\end{theorem}

$F(\omega ):=\{x\in C:x\in T(\omega ,x)\}\neq \emptyset .$

\textit{Then }$T$\textit{\ has a random fixed point theorem.}

\begin{proof}
Let us define $F:\Omega \rightarrow 2^{C}$ by $F(\omega )=\{x\in C:x\in
T(\omega ,x)\}.$ We will prove the measurability of $F.$ In order to do
this, we consider $B_{0}$ an arbitrary closed ball of $C,$ and let us denote

$L(B_{0}):=\tbigcap\limits_{k=1}^{\infty }\tbigcup\limits_{z\in B(B_{0};%
\frac{1}{k})}\{\omega \in \Omega :d(z,T(\omega ,z))<\frac{1}{k}\}.$

Now, we are proving that $F^{-1}(B_{0})=L(B_{0}).$

Firstly, let us consider $\omega \in F^{-1}(B_{0})$ and hence there exists $%
x_{0}\in B_{0}$ such that $x_{0}\in (T(\omega ,\cdot ))^{-1}(x_{0}).$

Since $T$ satisfies condition $\beta ,$ according to Lemma 8, for each $n\in
N,$ there exists $x_{n}\in B(B_{0},\frac{1}{n})$ such that $x_{n}\rightarrow
x_{0}$ when $n\rightarrow \infty $ and $d(x_{n},T(\omega ,x_{n}))<\frac{1}{n}%
.$ Therefore, $\omega \in L(B_{0})$ and then, $F^{-1}(B_{0})\subseteq
L(B_{0}).$

The rest of the proof is similar to the corresponding one in Theorem 3.1 in 
\cite{fie1}. Therefore, $F$ is measurable with nonempty closed values, and
according to the Kuratowski and Ryll-Nardzewski Proposition 1, $F$ has a
measurable selection $\xi :%
\Omega
\rightarrow C$ such that $\xi (\omega )\in T(\omega ,(\xi ,\omega ))$ for
each $\omega \in 
\Omega
$.
\end{proof}

\begin{corollary}
Let $(\Omega ,\mathcal{F})$ be a measurable space, $C$ be a closed separable
subset of a locally compact complete metric space and let $T:\Omega \times
C\rightarrow 2^{X}$ be a random operator which enjoys condition $\beta $ and
has closed values, such that for each $\omega \in \Omega ,$ $(T(\omega
.\cdot ))^{-1}:C\rightarrow 2^{C}$ is closed valued. Suppose, for each $%
\omega \in \Omega $ that the set
\end{corollary}

$F(\omega ):=\{x\in C:x\in T(\omega ,x)\}\neq \emptyset .$

\textit{Then }$T$\textit{\ has a random fixed point theorem.\medskip }

\begin{proof}
Since $T:\Omega \times C\rightarrow 2^{X}$ and for each $\omega \in \Omega ,$
$(T(\omega ,\cdot ))^{-1}:\Omega \times X\rightarrow 2^{C}$ are closed
valued and $C$ is locally compact, then $T$ fulfills condition ($\mathcal{P}$%
)$.$ We apply further Theorem 9.
\end{proof}

We will use further the following notation.

\begin{notation}
Let $X$ be a Fr\'{e}chet space, $C\subset X$ and $T:C\rightarrow 2^{X}.$ We
denote $S_{T}:C\rightarrow 2^{X}$ the correspondence defined by $S_{T}(x)=$co%
$\{T(x),x\}$ for each $x\in C.$
\end{notation}

We notice that for each $x\in X,$ $T(x)\subseteq S_{T}(x)$ and $%
T(x)=S_{T}(x) $ if only if $x\in T(x).$ We also notice that $x\in S_{T}(x)$
for each $x\in X.\medskip $

Now, we are defining the notion of local approximating pair of
correspondences.

\begin{definition}
Let $X$ be a topological vector space and $C\subset X$. Let $%
S,T:C\rightarrow 2^{X}$ be correspondences. We say that $(T,S)$ is a local
approximating pair if $T(x)\subseteq S(x)$ for each $x\in X$ and if $%
S(x_{0})=T(x_{0}),$ then, for each open neighborhood of the origin $V,$
there exists $U(x_{0})$ an open neighborhood of $x_{0}$ such that $%
S(x)\subseteq T(x)+V$ for each $x\in U(x_{0}).$
\end{definition}

Lemma 9 gives a new condition which implies (*).

\begin{lemma}
Let $C$ be a subset of topological vector space $X.$ Let $T:C\rightarrow
2^{X}$ be a correspondence such that $(T,S_{T})$ is a local approximating
pair. Then, $T$ satisfies (*).
\end{lemma}

\begin{proof}
Let us consider $x_{0}\in X$ such that $x_{0}\in T(x_{0})$, $\{x_{0}\}$
being included in a closed ball $B_{0},$ and the neighborhoods of the origin
of the following type: $V_{n}=B(0,\frac{1}{n}),$ $n\in \mathbb{N}.$ For each 
$n\in N,$ there exists a neighborhood $U_{n}$ of $x_{0}$ such that $%
U_{n}\subset B_{0}\cap B(x_{0},\frac{1}{n})$ such that $S_{T}(x)=$co$%
\{T(x)\cup x\}\subset T(x)+V_{n}$ for each $x\in U_{n}.$ For each $n\in 
\mathbb{N}$ we can pick $x_{n}\in U_{n}$ and then, for each $n\in \mathbb{N}$%
, $d(x_{n},x_{0})<\frac{1}{n}$ and $d(x_{n},T(x_{n}))<\frac{1}{n}.$
\end{proof}

Theorem 10 is easily obtained in the same way as Theorem 9.

\begin{theorem}
Let $(\Omega ,\mathcal{F})$ be a measurable space, $C$ be a closed separable
subset of Fr\'{e}chet space and let $T:\Omega \times C\rightarrow 2^{X}$ be
a random operator such that $(T,S_{T})$ is a local approximating pair$.$
Assume that $T$ enjoys condition ($\mathcal{P}$) and for each $\omega \in
\Omega $ the set $F(\omega ):=\{x\in C:x\in T(\omega ,x)\}\neq \emptyset .$
\end{theorem}

\textit{Then }$T$\textit{\ has a random fixed point theorem.\medskip }

\subsection{RANDOM FIXED POINT THEOREMS FOR UPPER SEMICONTINUOUS
CORRESPONDENCES}

The main aim of this subsection is to provide a new proof for the random Ky
Fan fixed point theorem. Our result is distinguished by the fact that $%
(\Omega ,\mathcal{F})$ is only a measurable space, without having other
additional properties. The proof is based on the upper approximation
technique, which is due to Ionescu-Tulcea \cite{ion}. The notions related to
this topic are presented below.\medskip

A correspondence $T:X\rightarrow 2^{Y}$ is \textit{quasi-regular (see }\cite%
{ion}) if: 1) it has open lower sections, that is $T^{-1}(y)$ is open in $X$
for each $y\in Y;$ 2) $T(x)$ is non-empty and convex for each $x\in X$ and
3) $\overline{T}(x)=$cl$T(x)$ for each $x\in X.$

The\textit{\ }correspondence $T$ is called \textit{regular} (see \cite{ion})
if it is quasi-regular and it has an open graph.\medskip

In \cite{ion} C. Ionescu Tulcea defined the notion of upper approximating
family for a correspondence.

Let $X$ be a non-empty set, $Y$ be a non-empty subset of topological vector
space $E$ and $T:X\rightarrow 2^{Y}.$ A family $(T_{j})_{j\in J}$ of
correspondences between $X$ and $Y$, indexed by a non-empty filtering set $J$
(denote by $\leq $ the order relation in $J$) is an \textit{upper} \textit{%
approximating family for }$T$ \cite{ion} if:

1) $T(x)\subset T_{j}(x)$ for each $x\in X$ and for each $j\in J;$

2) for each $j\in J$ there exists $j^{\ast }\in J$ such that, for each $%
h\geq j^{\ast }$ and $h\in J,$ $T_{h}(x)\subset T_{j}(x)$ for each $x\in X;$

3) for each $x\in X$ and $V\in $\ss , where \ss\ is a base of neighborhood
of 0 in $E$, there exists $j_{x,V}\in J$ such that $T_{h}(x)\subset T(x)+V$
if $h\in J$ and $j_{x,V}\leq h.\medskip $

By 1)-3) we deduce that:

4) for each $x\in X$ and $k\in J,$ $T(x)\subset \cap _{j\in J}T_{j}(x)=\cap
_{k\leq j,k\in J}T_{j}(x)\subset $ cl$T(x)\subset \overline{T}(x).\medskip $

Conditions for the existence of an approximating family for an upper
semicontinuous correspondence are given in the following Lemma. We recall
that a subset $X$ of a locally convex topological vector space\textit{\ }$E$
has \textit{the property} $(K)$ (see \cite{ion}) if, for every compact
subset $D$ of $X$, the set co$D$ is relatively compact in $E$.\medskip

\begin{lemma}
(see \cite{ion}). \textit{Let }$X$\textit{\ be a\ paracompact space and let }%
$Y$\textit{\ be a non-empty closed convex subset in a Hausdorff locally
convex topological vector space and with the property }$(K)$\textit{. Let }$%
T:X\rightarrow 2^{Y}$\textit{\ be a compact and upper semicontinuous
correspondence with non-empty convex compact values.\ Then, there exists a
filtering set }$J$\textit{\ such that there exists a family }$(T_{j})_{j\in
J}$\textit{\ of correspondences between }$X$\textit{\ and }$Y_{\text{ }}$%
\textit{\ with the following properties:}
\end{lemma}

1)\textit{\ for each }$j\in J,$\textit{\ }$T_{j}$\textit{\ is regular;}

2)\textit{\ }$(T_{j})_{j\in J}$\textit{\ and }$(\overline{T}_{j})_{j\in J}$%
\textit{\ are upper} \textit{approximating family for }$T;$

3)\textit{\ for each }$j\in J,\overline{T}_{j}$\textit{\ is continuous if }$%
Y_{\text{ }}$\textit{\ is compact.\medskip }

Theorem 11 is a random version of the Ky Fan fixed point Theorem. A new
proof is provided. Our result is distinguished by the fact that $(\Omega ,%
\mathcal{F})$ must be only a measurable space, without any additional
properties.

\begin{theorem}
Let $(\Omega ,\mathcal{F})$ be a measurable space, $C$ be a compact convex
separable subset of a Fr\'{e}chet space $X,$ and let $T:\Omega \times
C\rightarrow 2^{C}$ be a random operator. Suppose, for each $\omega \in
\Omega ,$ that $T(\omega .\cdot )$ is upper semicontinuous \textit{with
non-empty convex compact values}$.$
\end{theorem}

\textit{Then, }$T$\textit{\ has a random fixed point.}

\begin{proof}
Let us define $F:\Omega \rightarrow 2^{C}$ by $F(\omega )=\{x\in C:x\in
T(\omega ,x)\}.$ We will prove the measurability of $F.$ In order to do
this, we consider $B_{0}$ an arbitrary closed ball of $C,$ and let us denote

$L(B_{0}):=\tbigcap\limits_{k=1}^{\infty }\tbigcup\limits_{z\in B(B_{0},%
\frac{1}{k})}\{\omega \in \Omega :d(z,T(\omega ,z))<\frac{1}{k}\}.$

Now, we are proving that $F^{-1}(B_{0})=L(B_{0}).$

Firstly, let us consider $\omega \in F^{-1}(B_{0})$ and hence there exists $%
x_{0}\in B_{0}$ such that $x_{0}\in T(\omega ,x_{0}).$ Since $T(\omega
,\cdot )$ satisfies the assumptions of Lemma 10, there exists a filtering
set $J$\ such that there exists an upper approximating family $(T_{j}(\omega
,\cdot ))_{j\in J}$\ of correspondences between $C$\ and $X_{\text{ }}$\
such that for each $j\in J,$\ $T_{j}(\omega ,\cdot )$\ is regular. For $j\in
J,$ each $T_{j}(\omega ,\cdot )$ is lower semicontinuous, then it fulfils
condition (*)$.$ We will prove that condition (*) is also fulfilled by $%
T(\omega ,\cdot ).$

Since $(T_{j}(\omega ,\cdot ))_{j\in J\text{ }}$ is an approximating family
for $T(\omega ,\cdot ),$ for each $x\in C$ and $r>0$, there exists $%
j_{x,r}\in J$ such that $T_{i}(\omega ,x)\subset B(T(\omega ,x);r)$ if $i\in
J$ and $j_{x,r}\leq i.$ For $U=\tbigcup\limits_{k\geq 1}B(B_{0};\frac{1}{k}%
)\cap C$ and every $r>0,$ there exists $j_{U,r}\in J$ such that $%
T_{i}(\omega ,x)\subset B(T(\omega ,x);r)$ if $i\in J$ and $j_{U,r}\leq i.$

If $x_{0}\in T(\omega ,x_{0}),$ then, $x_{0}\in T_{i}(\omega ,x_{0})$ for
each $i\in J.$ For each $n\in \mathbb{N},$ there exists $j_{U,n}\in J$ and $%
i_{n}\in J$ such that $j_{U,n}\leq i_{n}$ and $T_{i_{n}}(\omega ,x)\subset
B(T(\omega ,x);\frac{1}{2n})$ for each $x\in U.$ Since $T_{i_{n}}(\omega
,\cdot )$ is lower semicontinuous, there exists $\{x_{n}^{i_{n}}\}$ such
that $d(x_{n}^{i_{n}},B_{0})<\frac{1}{2n}$ and $d(x_{n}^{i_{n}},T_{i_{n}}(%
\omega ,x_{n}^{i_{n}}))<\frac{1}{2n}.$ Consequently, $d(x_{n}^{i_{n}},T(%
\omega ,x_{n}^{i_{n}}))<\frac{1}{n}.$ Let us construct the sequence $%
\{z_{n}\}$ such that for each $n\in \mathbb{N},$ $z_{n}=x_{n}^{i_{n}}.$
Therefore, we found a convergent sequence $\{z_{n}\}$ such that $z_{n}\in
B(B_{0};\frac{1}{n})$ and $d(z_{n},T(\omega ,z_{n}))<\frac{1}{n}.$ Hence, $%
\omega \in L(B_{0})$ and the inclusion $F^{-1}(B_{0})\subseteq L(B_{0})$ is
proven$.$

Let $\omega \in \Omega ,$ let $B_{0}$ be a closed ball of $C$ and $\{z_{n}\}$
be a sequence in $C$ such that $d(z_{n},B)\rightarrow 0$ and $%
d(z_{n},T(\omega ,z_{n}))\rightarrow 0$ when $n\rightarrow \infty .$ The set 
$C$ is compact, then we can assume that $\{z_{n}\}$ is a convergent
sequence. Let $\lim_{n\rightarrow \infty }z_{n}=z_{0}\in B_{0}.$ Since $%
d(z_{n},T(\omega ,z_{n}))\rightarrow 0$ when $n\rightarrow \infty ,$ for
each $n\in \mathbb{N},$ there exists $y_{n}\in T(\omega ,z_{n})$ such that $%
d(z_{n},y_{n})\rightarrow 0$ when $n\rightarrow \infty .$ This fact assures
the convergence of the sequence $\{y_{n}\}$ and $\lim_{n\rightarrow \infty
}y_{n}=z_{0}.$ Now, we are using the upper semicontinuity of $T$ and we
conclude that for each $\varepsilon >0,$ there exists $N(\varepsilon )\in 
\mathbb{N}$ such that $T(\omega ,z_{n})\subseteq B(T(\omega
,z_{0});\varepsilon )$ for each $n>N(\varepsilon ).$ It follows that for
each $\varepsilon >0,$ $z_{0}\in B(T(\omega ,z_{0});\varepsilon )$ which
implies $z_{0}\in T(\omega ,z_{0}).$

We showed that the operator $T$ satisfies condition ($\mathcal{P}$).

We prove that $F:\Omega \rightarrow 2^{C}$ is measurable and has closed
values by following the same line as in the proof of Theorem 3.2 in \cite%
{fie1}.

In addition, we mention that since $C$ is compact convex and for each $%
\omega \in \Omega ,$ $T(\omega .\cdot )$ is upper semicontinuous with
non-empty convex compact values, then, according to Ky-Fan fixed point
theorem, the set $F(\omega ):=\{x\in C:x\in T(\omega ,x)\}\neq \emptyset .$

Consequently, $F$ is measurable with non-empty closed values. According to
the Kuratowski and Ryll-Nardzewski Proposition 1, $F$ has a measurable
selection $\xi :%
\Omega
\rightarrow C$ such that $\xi (\omega )\in T(\omega ,(\xi ,\omega ))$ for
each $\omega \in 
\Omega
$. We proved the existence of a random fixed point for $T.$
\end{proof}

We also establish a new random fixed point result concerning the upper
semicontinuous random operators.

\begin{corollary}
Let $(\Omega ,\mathcal{F})$ be a measurable space, $C$ be a compact convex
separable subset of a Fr\'{e}chet space $X,$ and let $T:\Omega \times
C\rightarrow 2^{C}$ be a random operator. Suppose, for each $\omega \in
\Omega ,$ that $(T(\omega .\cdot ))^{-1}:C\rightarrow 2^{C}$ is upper
semicontinuous \textit{with non-empty convex compact values}$.$
\end{corollary}

\textit{Then, }$T$\textit{\ has a random fixed point.}

\begin{proof}
According to Theorem 2, there exists a measurable mapping $\xi :\Omega
\rightarrow C$ such that\textit{\ }for every $\omega \in \Omega $%
\c{}
$\xi (\omega )\in (T(\omega ,\cdot ))^{-1}(\xi (\omega ),$ that is, for
every $\omega \in \Omega $%
\c{}
$\xi (\omega )\in T(\omega ,\xi (\omega )).$ Therefore, we obtained a random
fixed point for $T.$
\end{proof}

In the end of the paper we state the following theorem, which is a
consequence of Theorem 11.

\begin{notation}
Let $\Omega $ be a nonempty set, $X,Y$ be topological spaces, $C$ be a
closed subset of $X$ and let $T:\Omega \times C\rightarrow 2^{X}$ be an
operator. We denote $\overline{T}:\Omega \times C\rightarrow 2^{C}$ the
operator defined in the following way: for each $\omega \in \Omega $ and $%
x\in C,$ $\overline{T}(\omega ,x):=\{y\in C:(x,y)\in $cl$_{C\times C}$ Gr$%
T(\omega ,\cdot )\}$
\end{notation}

\begin{theorem}
Let $(\Omega ,\mathcal{F})$ be a measurable space, $C$ be a compact convex
separable subset of a Fr\'{e}chet space $X,$ and let $T:\Omega \times
C\rightarrow 2^{C}$. Suppose that $\overline{T}:\Omega \times C\rightarrow
2^{C}$ is a random operator and for each $\omega \in \Omega ,$ $\overline{T}%
(\omega .\cdot )$ has\textit{\ non-empty convex values, where }$\overline{T}%
(\omega .\cdot )$ is defined by $\overline{T}(\omega ,x):=\{y\in C:(x,y)\in $%
cl$_{C\times C}$Gr$T(\omega ,\cdot )\}.$
\end{theorem}

\textit{Then, there exists a measurable mapping }$\xi :%
\Omega
\rightarrow C$ \textit{such that} $\xi (\omega )\in \overline{T}(\omega
,(\xi ,\omega ))$ for each $\omega \in 
\Omega
$\textit{.}

\begin{proof}
For each $\omega \in \Omega ,$ $\overline{T}(\omega .\cdot )$ is upper
semicontinuous with non-empty convex compact values, then, we can apply
Theorem 11.
\end{proof}

\section{CONCLUDING\ REMARKS}

We have proven the existence of random fixed points for almost lower
semicontinuous and lower semicontinuous operators defined on Fr\'{e}chet
spaces. Our research extends on some results which exist in literature. It
is an interesting problem to find new types of operators which satisfy weak
continuity properties and have random fixed points.

\end{document}